\documentclass[hidelinks,onefignum,onetabnum]{siamart250211}

\usepackage{caption}
\usepackage{lipsum}
\usepackage{amsfonts}
\usepackage{graphicx}
\usepackage{epstopdf}
\usepackage{geometry}
\usepackage{amsmath,amssymb}
\usepackage{algorithm}
\usepackage{algorithmic}
\usepackage{bm}
\usepackage[caption=false]{subfig}
\usepackage{appendix}
\usepackage{multirow}
\usepackage{braket}
\usepackage[english]{babel}
\usepackage{hyperref}
\usepackage[capitalize]{cleveref}
\usepackage[square,numbers,sort&compress]{natbib}
\usepackage{adjustbox}
\usepackage{xspace}
\ifpdf
  \DeclareGraphicsExtensions{.eps,.pdf,.png,.jpg}
\else
  \DeclareGraphicsExtensions{.eps}
\fi

\def\le{\leqslant}
\def\ge{\geqslant}

\numberwithin{equation}{section}

\newsiamremark{remark}{Remark}
\newsiamremark{hypothesis}{Hypothesis}
\crefname{hypothesis}{Hypothesis}{Hypotheses}
\newsiamthm{claim}{Claim}
\newsiamthm{example}{Example}

\headers{Uniform Semiclassical Observable Error Bound of Trotterization}{Di Fang and Conrad Qu}

\title{
Uniform semiclassical observable error bound of Trotter-Suzuki splitting: a simple algebraic proof
\thanks{
\funding{This work is supported by National Science Foundation via the grant DMS-2347791 and DMS-2438074 (D.F. and C.Q.), and the U.S. Department of Energy, Office of Science, Accelerated Research in Quantum Computing Centers, Quantum Utility through Advanced Computational Quantum Algorithms, grant no. DE-SC0025572 (D.F.). }}}

\author{Di Fang\thanks{Department of Mathematics and Duke Quantum Center, Duke University 
  (\email{di.fang@duke.edu}).}
\and Conrad Qu \thanks{Department of Mathematics and Department of Electrical and Computer Engineering, Duke University
  (\email{conrad.qu@duke.edu}).}}

\usepackage{amsopn}

\newcommand{\REV}[1]{{#1}}

\newtheorem{hyp}{Assumption}

\usepackage{url}
\usepackage{enumitem}
\usepackage{indentfirst}

\numberwithin{equation}{section}

\newcommand{\I}{\mathrm{i}}

\newcommand{\norm}[1]{\left\lVert#1\right\rVert}

\newcommand{\Or}{\mathcal{O}}

\ifpdf
\hypersetup{
  pdftitle={
  Uniform semiclassical observable error bound of Trotter-Suzuki splitting: a simple algebraic proof
  },
  pdfauthor={Di Fang, Conrad Qu}
}
\fi

\crefname{section}{Section}{Sections}
\Crefname{section}{Section}{Sections}
\crefname{appendix}{Appendix}{Appendices}
\Crefname{appendix}{Appendix}{Appendices}
\crefname{figure}{Figure}{Figures}
\Crefname{figure}{Figure}{Figures}
\crefname{table}{Table}{Tables}
\Crefname{table}{Table}{Tables}
\crefname{theorem}{Theorem}{Theorems}
\Crefname{theorem}{Theorem}{Theorems}
\crefname{lemma}{Lemma}{Lemmas}
\Crefname{lemma}{Lemma}{Lemmas}
\crefname{proposition}{Proposition}{Propositions}
\Crefname{proposition}{Proposition}{Propositions}
\crefname{corollary}{Corollary}{Corollaries}
\Crefname{corollary}{Corollary}{Corollaries}

\begin{document}

\maketitle

\begin{abstract}
Efficient simulation of the semiclassical Schrödinger equation has garnered significant attention in the numerical analysis community. While controlling the error in the unitary evolution or the wavefunction typically requires the time step size to shrink as the semiclassical parameter $h$ decreases, it has been observed--and proved for first- and second-order Trotterization schemes--that the error in certain classes of observables admits a time step size independent of $h$. In this work, we explicitly characterize this class of observables and present a new, simple algebraic proof of uniform-in-$h$ error bounds for arbitrarily high-order Trotterization schemes. Our proof relies solely on the algebraic structure of the underlying operators in both the continuous and discrete settings. Unlike previous analyses, it avoids Egorov-type theorems and bypasses heavy semiclassical machinery. To our knowledge, this is the first proof of uniform-in-$h$ observable error bounds for Trotterization in the semiclassical regime that relies only on algebraic structure, without invoking the semiclassical limit.
\end{abstract}
 


\section{Introduction}
Simulation of quantum dynamics has been a foundational motivation for the development of quantum computers~\cite{Feynman1982}, and remains one of the most promising applications. In particular, quantum algorithms for Hamiltonian simulation aim to efficiently approximate the time evolution of quantum systems, and are widely regarded as a core primitive in quantum algorithm design. However, as dictated by the no-fast-forwarding theorem, e.g., \cite[Theorem 3]{BerryAhokasCleveEtAl2007}, \cite[Theorem 5]{Childs2010CMP}, \cite{Kothari2010}, the cost of simulation in the worst case is expected to scale at least linearly with the time and norm of the Hamiltonian, implying a greater challenge in addressing multiscale quantum systems where certain small parameters induce separated scales and cause the Hamiltonian norm to grow as they vary.

A prominent example of such scale separation behavior is the so-called semiclassical Schrödinger equation:
\begin{equation}\label{eq:semiclassical_schd}
    i h \partial_t u^h(t,x) = -\frac{h^2}{2} \Delta u^h(t,x) + V(x) u^h(t,x),
\end{equation}
where $h >0 $ is the (dimensionless) reduced Planck constant, $x \in \mathbb{R}^d$ for dimension $d$, and $V(x)$ is the potential. In many contexts--especially in molecular quantum dynamics--this equation arises in the semiclassical regime, where $h \ll 1$ denotes the semiclassical parameter that quantifies the underlying scale separation. For example, in molecular dynamics, $h$ represents the square root of the mass ratio between electrons and nuclei. This stands in contrast to quantum dynamics of electrons (without scale separation), where $h$ is typically set to one in atomic units.

This equation emerges naturally in molecular quantum dynamics, particularly under the Born-Oppenheimer approximation (see its derivation in, e.g.,~\cite{LasserLubich2020,Lubich2008,Hagedorn2007}).
Beyond its origin in chemical and physical applications, the semiclassical Schrödinger equation has recently been found to be useful in designing efficient quantum algorithms for optimization~\cite{ZhangLengLi2021,LiuSuLi2023,LengHickmanLiWu2023}.

From a computational standpoint, simulating~\cref{eq:semiclassical_schd} is particularly challenging due to the oscillatory nature introduced by the small parameter $h$. In order to resolve accurately the wavefunction, both the spatial grid size $\Delta x$ and the time step size $\Delta t$ must decrease as $h \to 0$, significantly increasing computational cost. This difficulty has led to a large body of work in the numerical analysis community, focusing on the development of more efficient trajectory-based methods that leverage the rich mathematical structure of the $h \to 0$ limit (see, e.g., reviews \cite{LasserLubich2020,JinMarkowichSparber2011}). Recently, however, grid-based methods in the semiclassical regime have regained attention, driven by the advent of efficient quantum algorithms for Hamiltonian simulation~\cite{BornsWeilFang2022,JinLiLiu2021}. While the need for small $\Delta x$ increases the dimension of the discretized Hilbert space, quantum computers can handle this efficiently,
resulting in only $\mathrm{polylog}(h^{-1})$ cost overhead. However, the shrinking time step $\Delta t$ still contributes to the quantum simulation cost polynomially, leading to a total cost scaling as $\mathrm{poly}(h^{-1})$, which is consistent with the intuition behind the no-fast-forwarding theorem.

Despite this challenge of efficient time stepping, it was first observed in~\cite{BaoJinMarkowich2002} that when one is interested in computing (quadratic) physical observables (as opposed to full wavefunctions), time-splitting spectral methods with large time step size $\Or(1)$ can still yield surprisingly accurate results. This striking phenomenon has sparked a growing body of work investigating efficient algorithms for observable accuracy in various contexts (see, e.g., ~\cite{BaoJinMarkowich2003,BaoJakschMarkowich2003,Carles2013,CarlesGallo2017,FaouLubich2004,JinWuYang2008,FaouGradinaruLubich2009,BaderIserlesKropielnickaSingh2014,FangJinSparber2018,CaiFangLu2022,FilbetGolse2024,FilbertGolse2025} and those cited by the reviews~\cite{JinMarkowichSparber2011,BaoCai2013,LasserLubich2020}).
Other significant developments include advances in the observability theory of semiclassical Schrödinger dynamics~\cite{GolsePaul2021}, further motivating the design of simulation methods tailored to observable quantities.

In recent years, significant progress has been made toward understanding the observable accuracy of Trotter-based methods in semiclassical quantum dynamics. A central challenge in this setting is to establish error bounds that remain uniform with respect to the semiclassical parameter $h$, while retaining the nominal convergence rate of the splitting scheme. Let $\Delta t$ denote the Trotter step size, so that the number of steps to simulate time $t$ is approximately $t/\Delta t$.

One significant line of work analyzes second-order Trotter formulas for the von Neumann equation by comparing the Husimi transforms of the exact and approximate quantum states using the quantum Wasserstein distance~\cite{GolseJin19}. This approach yields an observable error bound of $\Or(\Delta t^2 + h^{1/2})$, and under additional assumptions on the initial state, a uniform-in-$h$ bound of $\Or(\Delta t^{2/3})$. For the semiclassical Schrödinger equation, a rigorous error estimate on observables using a second-order splitting method was established in~\cite{LasserLubich2020}, a foundational work that has since inspired numerous follow-up studies. Leveraging Egorov theorems for Husimi functions and the symplectic structure of the St\"ormer-Verlet integrators, it shows an error bound of $\Or(\Delta t^2 + h^2)$ when the observables are quantizations of Schwartz functions. Further advance was made in~\cite{FangTres2021}, where uniform bounds were extended to weakly nonlinear quantum-classical models, improving the rate to $\Or(\Delta t^{4/3})$ under structural assumptions on the initial wavepacket. 

Notably, these analyses have largely focused on spatially continuous cases and rely rather heavily on semiclassical tools such as Egorov-type theorems, Husimi functions, and Wigner transforms. Moreover, due to additive scaling in $h$, these results are most relevant in low-precision regimes where the desired accuracy exceeds $\Or(t h^2)$, as discussed in~\cite{JinLiLiu2021}. In this overview, we restrict attention to linear systems and to improvements in observable error bounds—precisely the setting of our work. We do not further review results for nonlinear systems, but we note that, under WKB‐type initial data, uniform-in-$h$ estimates for position and current densities have been obtained for several nonlinear models~\cite{Carles2013,CarlesGallo2017}.

At a conceptual level, it is also worth noting that additive error bounds -- such as those scaling with $h^2$ or $h^{1/2}$ -- may be somewhat unphysical in this context. These bounds suggest that the error persists even as the number of Trotter steps tends to infinity, indicating a deviation that stems not from the numerical approximation but from the analytical framework used to interpret it. In particular, such additive-in-$h$ bounds often arise from comparing the quantum observable to its semiclassical limit, even though the simulation itself does not take that limit. As a result, such bounds may fail to reflect the genuinely quantum characteristics of the system. Specifically, if the error depends additively on $h$, then simulating the corresponding classical dynamics (i.e., the semiclassical limit) would produce similar accuracy as the Trotterization for the full quantum dynamics. This suggests that the quantum simulation may not be meaningfully capturing the quantum nature of the system.

More recently, \cite{BornsWeilFang2022} provided the first observable error bound for the semiclassical Schrödinger equation that is both uniform in $h$ and preserves the formal convergence rate of the splitting method. Crucially, the analysis was further carried out in the spatially discretized setting, making it applicable to numerical simulations in practice. It also extends the class of admissible observables from the Schwartz class, as considered in prior state-of-the-art results, to a broader symbol class. However, that result was limited to first- and second-order Trotter schemes and relied fundamentally on advanced techniques from semiclassical and discrete microlocal analysis. These methods, while powerful, are difficult to generalize to higher-order schemes, highlighting the need for a more elementary and extensible approach.

In summary, although prior results have provided significant insights into the observable accuracy of low-order Trotter-based methods, new and simpler proofs remain highly desirable to deepen theoretical understanding and broaden applicability. This leads to the open question:

\begin{center}
\textit{Can we rigorously establish a uniform-in-$h$ observable error bound for Trotterization--at arbitrary high order--using a simple algebraic framework, without invoking Egorov-type theorems?}
\end{center}

We provide an affirmative answer by presenting a new error analysis based solely on Taylor expansion of observables and direct commutator calculations in the underlying Lie algebra of operators. Our approach applies both to continuous and discretized systems and yields a clean algebraic proof of the observable error bound, achieving uniformity in the semiclassical parameter $h$ without sacrificing the high-order accuracy of the splitting scheme.

\medskip

\textbf{Our Contributions:}

We establish a uniform-in-$h$ observable error bound for arbitrary high-order Trotterization, which, to the best of our knowledge, is new in the literature. In terms of methodology, we leverage key algebraic structures of the underlying operators--namely, the height-reduction and width-expansion properties--which may be of independent interest. Our proof circumvents the need for semiclassical microlocal analysis -- such as Egorov-type theorems -- and instead relies entirely on elementary algebraic arguments, including Taylor expansions of observables and direct commutator calculations in the Lie algebra. This algebraic approach offers a more accessible and broadly applicable framework for analyzing Trotterized semiclassical simulations and opens the door to extensions beyond the limitations of existing Egorov-based analyses. The resulting  proof is remarkably simple yet elegant, offering several key advantages: it greatly simplifies the analysis of high-order Trotter schemes, extends naturally to time-dependent potentials $V(t, x)$, and applies uniformly to both spatially continuous operators and fully discretized settings through explicit algebraic arguments.
\begin{figure}[htp]
    \centering
    {\includegraphics[width=.6\textwidth]{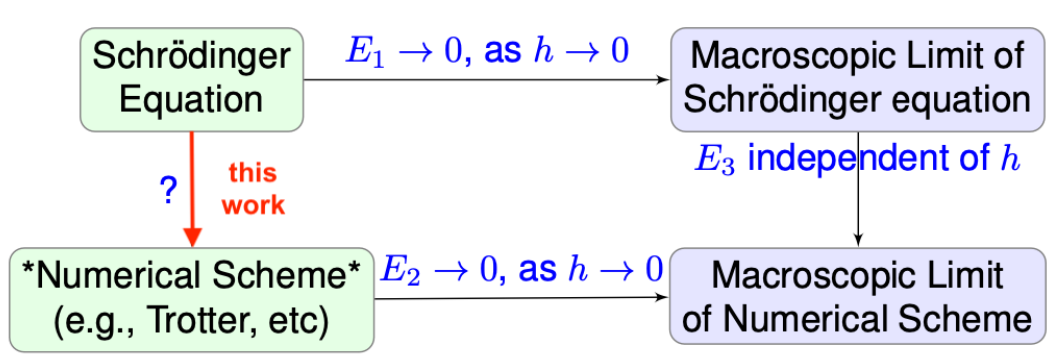}}
    \label{fig:commutators}
    \caption{This figure illustrates the distinct proof strategy employed in our work compared to prior significant advancements~\cite{BaoJinMarkowich2002,JinMarkowichSparber2011,GolseJin19,LasserLubich2020,FangTres2021,JinLiLiu2021} on this topic. Rather than using macroscopic (semiclassical) limits as an intermediate step and controlling the total error via a decomposition into multiple contributions, we directly estimate the error in the time-evolved observable produced by the numerical scheme. $E_1$, $E_2$, and $E_3$ denote the errors between the two objects they connect, measured in an appropriate error metric (see the cited work above for further details). Different from the previous analysis developed by one of the authors for first- and second-order Trotterizations, this work develops a new purely algebraic proof and extends the uniform-in-$h$ error bounds to arbitrarily high-order Trotterization schemes.}
\end{figure}
On an informal level, we have showed that the time steps of Trotterization can be made independent of $h^{-1}$.  In particular, we proved a uniform-in-$h$ observable error bound of \REV{$p$}-th order Trotter formula and polynomial observables $O$, i.e., (informally; see~\cref{thm:global_error} for detailed definition)
\[
 \norm{ \mathcal{U}_\text{app}^\dagger O \mathcal{U}_\text{app} - e^{i H t} O e^{-\I H t}} \leq C n^{-p},
 \]
where $C$ is \textit{independent} of $h^{-1}$ (!), and $n$ is the number of Trotter steps. 

\smallskip
\textbf{Organization:}

The rest of the paper is organized as follows.  In \cref{sec:problem_setup} we introduce our notation, recall the semiclassical Schrödinger operator and the $p$-th–order Trotter–Suzuki decomposition, and state our two main results: \REV{a} uniform-in-$h$ observable-error bound for arbitrary high-order Trotter formulas and a nested-commutator error estimate valid for general Hamiltonians $H=A+B$ and observables $O$.  \cref{sec:continuous_trotter} gives an algebraic Taylor-expansion proof of the local truncation error of the $p$-th order Trotter formula,
with full control of the $h$-dependence.  In \cref{sec:algebraic_structures}, we analyze the underlying algebraic structures: we prove height-width lemmas in the semi-discrete case, relevant to spectral discretizations, and extend it to finite differences to obtain uniform-in-$h$ commutator bounds.  \cref{sec:numerics} presents numerical results confirming our theoretical findings.  We conclude in \cref{sec:conclusion} with a discussion of implications and future directions.

\section{Problem Setup and Main Result} 
\label{sec:problem_setup}
In this section, we set up the problem by reviewing the semiclassical Schrödinger operator and Trotterization, and then present the main uniform estimates we obtained.

\subsection{Higher-Order Trotter Formula and Semiclassical Notations}

We revisit the $p$-th order Trotter-Suzuki decomposition applied to simulate the time evolution operator  $U(t) = e^{-iHt}$ for a Hamiltonian $H = A + B$, where $A$ and $B$ are its distinct Hamiltonian components.
\REV{In the work, all Hamiltonians $A$, $B$, and $H$ are Hermitian operators, so $e^{-iAt}$, $e^{-iBt}$, and $e^{-iHt}$ are unitary.}

The general expression for the $p$-th order Trotter-Suzuki formula ($p = 2k$ is an even integer)~\cite{Suzuki1993}, denoted as $ U_p^{\Delta t} =  U_p(\Delta t) $, is given recursively. For the 2nd-order case, the Trotter-Suzuki formula simplifies to:
\begin{equation}\label{eq:suzuki_2}
    U_2^{\Delta t} = e^{-i \frac{\Delta t}{2} A} e^{-i \Delta t B} e^{-i \frac{\Delta t}{2} A}.
\end{equation}
For higher even orders, specifically the $2k$-th order, this structure is extended recursively as:

\begin{equation}\label{eq:suzuki_2k}
    U_{2k}^{\Delta t} := \REV{\bigl(U_{2k-2}^{\,u_k \Delta t}\bigr)^{2}\, U_{2k-2}^{\,(1-4u_k)\Delta t}\, \bigl(U_{2k-2}^{\,u_k \Delta t}\bigr)^{2}},
\end{equation}
where $u_k = \frac{1}{4 - 4^{1/(2k-1)}}$. 
Note that the \REV{$2k$}-th order Suzuki construction contains $2 \cdot 5^{k-1} + 1$ many exponentials. Alternatively, the Yoshida construction~\cite{Yoshida1990} offers a comparable formulation for the \REV{$2k$}-th order that contains $2 \cdot 3^{k-1} + 1$ terms. The number of exponentials (stages) in \cref{eq:suzuki_2k} is denoted as $l$. Note that $p=2$ is also called the Strang splitting in numerical analysis literature. We also remark that the Trotter–Suzuki formula fits within the broader framework of composition methods; see, e.g., the book~\cite{BlanesCasas2025} and the review~\cite{BlanesCasasMurua2008}, together with the references cited therein, for further background.

In the following, unless otherwise specified, we consider the semiclassical Schrödinger equation~\cref{eq:semiclassical_schd} for $0 < h \leq 1$ with periodic boundary conditions. For notational simplicity, we define
\begin{equation}\label{eq:semiclassical_A_B}
A = -\frac{h}{2} \Delta, \quad B = \frac{1}{h}V(x)
\end{equation}
or their discretizations. We note that interchanging the order of $A$ and $B$ in our setting does not affect the uniform-in-$h$ result, thanks to the underlying algebraic structure and the regularity of $V$. We emphasize that although the explicit dependence on $h$ is omitted in the notation for $A$ and $B$, both operators--and their spatial discretizations--depend on $h$. Throughout the paper, we assume that the real-valued potential $V$ satisfies the following assumption:
\begin{hyp}\label{assump:V}
$V$ is smooth and bounded together with all its derivatives, i.e. $V \in S(1)$.
\end{hyp}
Note that the notation $S(1)$ is commonly used in semiclassical analysis~\cite{zworski_book}, and adopt it here for notational simplicity.
The constant $C$ in \cref{thm:global_error} depends on the derivatives of $V$, stemming from the algebraic calculations in~\cref{sec:algebraic_structures}.
As our goal is to discuss arbitrarily high-order Trotter–Suzuki splittings, we impose this regularity assumption on $V$ for simplicity. Relaxing this assumption on $V$ is an interesting direction for future work. From a numerical perspective, one essentially needs $V$ to be bounded together with its derivatives on the computational domain. For example, for the harmonic oscillator potential $V = x^2$ with strong confinement, the effective support of a wavefunction is finite, so the potential and its derivatives are effectively bounded on the computational domain. For such a potential, the operators in \cref{eq:semiclassical_A_B}, as well as their sum, are self-adjoint (from their respective domains to $L^2$), and so are their discretized versions (using either finite difference or spectral methods) with periodic boundary conditions. Consequently, all exponentials appearing in the Trotter–Suzuki formulas~\cref{eq:suzuki_2,eq:suzuki_2k} are unitary.

\subsection{Main Results}

We define the exact evolution of an observable $O$ under the Hamiltonian $H$ as:
\begin{equation}\label{eq:T_exact_def}
T(t) := e^{iHt} O e^{-iHt}.
\end{equation}
The time evolution of $ O $ under the $p$-th order Trotter formula, after $ n $ time steps of size $ \Delta t $, is given by:
\begin{equation}\label{eq:T_pn_def}
T_{p,n}(\Delta t) := ({U_p^{\Delta t}}^\dagger)^n O (U_p^{\Delta t})^n,
\end{equation}
where $U_p^{\Delta t} $ is the $p$-th order Trotter approximation to $ e^{-iH\Delta t}$. Here $^\dagger$ denotes the Hermitian conjugate, following the notational convention in quantum mechanics. For matrices, this is the conjugate transpose; for continuous operators, it is the adjoint with respect to the $L^2$ inner product.

For the semiclassical Schr\"odinger setting considered here, we consider the set of polynomial observables given by
\begin{hyp}\label{assump:observable}
We consider observables of the form
\begin{equation}\label{eq:poly_O_form}
    O = \sum_{m=0}^q y_m(x)\, h^m \partial_x^m, \quad y_m \in S(1)\ \text{for all } m.
\end{equation}
\end{hyp}
This essentially implies that, after suitable discretization, the observable has an operator norm of order $h^0$; see the discussion preceding \cref{thm:nested_poly_obs}. As is standard in numerical simulations of \cref{eq:semiclassical_schd}, we consider spatial discretizations with periodic boundary conditions. In this case, we assume the functions $y_m$ are periodic as well, ensuring that their derivatives remain bounded at the boundaries.

We now state the central theorem of bounding the global error of the $p$-th order Trotter formula. The proof of this theorem follows from the analysis of the local error in \cref{sec:continuous_trotter} and \cref{sec:algebraic_structures}.

\begin{theorem}[Global observable error for $p$-th order Trotterization]
\label{thm:global_error}
Let $A$ and $B$ be the discrete operators obtained from the finite-difference (or spectral) discretization of the semiclassical Schr\"odinger operators as in~\cref{eq:semiclassical_A_B} for $x \in \mathbb{R}$, with the real-valued potential $V$ satisfying Assumption~\ref{assump:V}. Let $H = A + B$ be the resulting discrete Hamiltonian.
Let $O$ be the discrete observable obtained by discretizing the observables specified in Assumption~\ref{assump:observable} using the same numerical scheme.
Let $T(t)$ and $T_{p,n}(\Delta t)$ be the time-evolved observables under the exact dynamics and the Trotter--Suzuki splitting, respectively, in the Heisenberg picture, as defined in~\cref{eq:T_exact_def,eq:T_pn_def}. Let $t = n \Delta t = \Or(h^{0})$. Then the global time-evolved observable error in operator norm satisfies
\begin{equation}
    \|\, T_{p,n}(\Delta t) - T(t) \| \;\le\; C\, \Delta t^p,
\end{equation}
where $\|\cdot\|$ denotes the spectral norm, and $C$ is a constant independent of $h^{-1}$ while depending on $t$, the observable $O$ and the potential $V$ (through derivatives), and the Trotter-Suzuki order $p$.
\end{theorem}

\cref{thm:global_error} follows from the local error bound, which we will derive in 
\cref{sec:continuous_trotter} using a Taylor series expansion of the exact and Trotter-evolved observables. \cref{sec:algebraic_structures} further establishes that the commutator scaling is independent of the semiclassical parameter $h$.

\begin{remark}[On the Dimensionality]
In principle, the result holds for general spatial dimension $d$. 
However, for notational simplicity, the algebraic structure underlying our proof is developed only in one spatial dimension. 
The extension to higher dimensions is straightforward: since 
$
\Delta = \sum_{j=1}^d \partial_{x_j}^2,
$
the same algebraic argument applies to each component $x_j$ independently. 
Correspondingly, in the height definition in \cref{eq:def_ht_wd}, the derivative $\partial_x^m$ should be interpreted as a multi-index derivative with multi-index $m$, and its height becomes $|m|$, the sum of its components. 
For mathematical rigor, we therefore state the main theorem only in one spatial dimension, while the proof extends verbatim to higher dimensions.
\end{remark}

An immediate consequence of Theorem~\ref{thm:global_error} is an estimate for the number of Trotter steps required to achieve a given precision $\epsilon$ in the operator norm. This leads to the following corollary.

\begin{corollary}[Query Complexity of Trotter-Suzuki Splittings]\label{cor:query_complexity}
Under the same conditions as in Theorem~\ref{thm:global_error}, let $n_p$ denote the total number of time steps required for the $p$-th order Trotter-Suzuki splitting to achieve an $\epsilon$-approximation for the time-evolved observable. Then for $t = \Or(h^0)$, we have
$
    n_p = \REV{\Or}\!\left(h^0 \,\epsilon^{-1/p}\right)$, 
that is, $n_p = \Or(\epsilon^{-1/p})$
independent of the semiclassical parameter $h$.
\end{corollary}
\begin{remark}
    On the other hand, if one evaluates the error in the wave function or in the unitary evolutions, the number of Trotter steps needed has a polynomial overhead in $h^{-1}$ (see, e.g., \REV{the insightful work \cite{DescombesThalhammer2010} for a rigorous proof for $p$-th order splitting, as well as} \cite{JinLiLiu2021}; see also \cite[Table~1 and Appendix~A]{BornsWeilFang2022} for a discussion of the $\operatorname{poly}(h^{-1})$ query dependence in other quantum algorithms). 
\end{remark}

As a byproduct of the proof of our main theorem (\cref{thm:global_error}), we obtain the following general bound on observable expectation errors under Trotterization: any $p$-th order Trotter formula incurs an observable error that can be bounded by $(p+1)$ layers of nested commutators involving the summands of the Hamiltonian and the observable $O$. We stress that, unlike the main theorem, this result applies to general Hamiltonians and observables, beyond the semiclassical regime. In quantum simulation applications, Trotterization is typically analyzed in a setting where all Hamiltonians are bounded operators acting on the same Hilbert space -- most commonly finite-dimensional matrices -- and each component (such as $A$ and $B$) is bounded as well.  
Any quantum state is represented by the ket-notation $\ket{\psi}$, which is always normalized so that $\|\ket{\psi}\|=1$ in the 2-norm.  
Under this normalization, the observable expectation value of any bounded observable $X$ obey
\[
|\bra{ \psi} X \ket{\psi}| \leq \norm{\ket{\psi}} \norm{X \ket{\psi}} \leq \|X\| \norm{\ket{\psi}}^2 = \norm{X}.
\]
This allows the observable expectation error to be controlled directly by the operator-norm error.
This yields the following useful result, which may be of independent interest in quantum computing.
To make it accessible to that community, we state the next proposition using the notational convention standard in quantum algorithms, even though it is not fully consistent with the notation used elsewhere in the paper. Quantum-leaning readers can read this result independently of the semiclassical sections, while numerical-analysis readers may safely skip it.
\begin{proposition}[General observable error bound for $p$-th order Trotter-Suzuki formula]\label{prop:general_observable_error}
Let $t>0$ and let $\ket{\psi_0}$ be a quantum state.
For any bounded observable $O$ and bounded Hamiltonian $H = A + B$, the Heisenberg-picture observable satisfies
\begin{align}
\Bigl|
    \bra{\psi_0}
        U_p^{t}{}^\dagger \, O \, U_p^{t}
    \ket{\psi_0}
    -
    \bra{\psi_0}
        e^{iHt} O e^{-iHt}
    \ket{\psi_0}
\Bigr|
\le \|U_p^{t}{}^\dagger \, O \, U_p^{t} - e^{iHt} O e^{-iHt} \|
\le C_p \, \beta_{\mathrm{comm}} \, t^{p+1},
\end{align}
where $U_p^t$ is the $p$-th order Trotter-Suzuki formula defined in \cref{eq:suzuki_2,eq:suzuki_2k},
$C_p$ depends only on the Trotter order $p$ (and the Suzuki coefficients), and
\begin{equation}\label{eq:beta}
\beta_{\mathrm{comm}}
:=
\max_{\substack{
q_1,\dots,q_k:\; q_1+\cdots+q_k = p+1 \\
H_j \in \{A,B\},\; j=1,\dots,k
}}
\bigl\|
    \mathrm{ad}_{H_k}^{q_k}
    \cdots
    \mathrm{ad}_{H_1}^{q_1}(O)
\bigr\|.
\end{equation}
\end{proposition}

\section{Trotter Error for Time-Evolved Observables}
\label{sec:continuous_trotter}
In this section, we establish the local error bound for the $p$-th order Trotter formula when applied to time-evolved observables. The central result of this section is the following theorem, which provides a uniform-in-$h$ upper bound on the local error for the semiclassical Schr\"odinger operators.

\begin{theorem}[Local error of observable for $p$-th order Trotter]\label{thm:trot_p} 
Let $ O $ be an observable of order $ h^0 $ as defined in \cref{eq:poly_O_form} by finite difference (or spectral) spatial discretizations. Denote the exact time evolution of the observable under the unitary $T(\Delta t)$ as given in \cref{eq:T_exact_def}, and the one-step time evolution under the $p$-th order Trotter formula is given by
\begin{equation}\label{eq:ob_trotter_p}
    T_p(\Delta t) := {U_p^{\Delta t}}^\dagger O U_p^{\Delta t},
\end{equation}
where $ U_p^{\Delta t} $ is the $p$-th order Trotter approximation of the time evolution operator.
Then the difference between the exact and approximate evolutions of the observable can be estimated as:
\begin{equation}\label{eq:ob_err_trotter_p}
    \norm{ T_p(\Delta t) - T(\Delta t)}\leq C \Delta t^{p+1},
\end{equation}
where $C$ is a constant independent of $h^{-1}$, and $\Delta t$ is the time step size.
\end{theorem}

In the following proof, we begin by analyzing the error for a general bounded Hamiltonian $H = A + B$ and a bounded observable $O$, which will later correspond to discretized matrices in the semiclassical setting. In \cref{subsec:comm_ob}, we establish that the observable error exhibits commutator scaling, as stated in \cref{prop:general_observable_error}. This algebraic structure
is then analyzed in detail in \cref{sec:algebraic_structures}, both in the spatially continuous setting--relevant for spectral methods--and in the spatially discrete setting for finite difference schemes, ultimately yielding the desired estimate.

\subsection{Taylor Series Expansion of Exact Evolution}
We begin by expanding the exact time-evolved observable $O$ under the Hamiltonian $H = A + B$ with bounded operators $O$, $A$, and $B$. Recall the observable at time $s$ is given by:
\begin{equation}
T(s) = e^{i H s} O e^{-i H s}.
\end{equation}
This operator can be expressed as a Taylor series expansion around  $s = 0$. For clarity, the action of $H$ on $O$ is represented in terms of the adjoint operator $\mathrm{ad}_H(O)$, defined recursively as:
\begin{equation} \label{eq:ad_def}
\mathrm{ad}_H(O) := [H, O], \quad \mathrm{ad}_H^n(O) := [H, \mathrm{ad}_H^{n-1}(O)], \quad n\geq 2.
\end{equation}
Using this definition, the expansion for $T(s)$ takes the form:
\[
T(s) = O + i s \mathrm{ad}_H(O) + \frac{(i s)^2}{2!} \mathrm{ad}_H^2(O) + \cdots + \frac{(i s)^p}{p!} \mathrm{ad}_H^p(O) + \widetilde{\mathcal{R}}_p(s),
\]
where $\widetilde{\mathcal{R}}_p(s)$ denotes the remainder terms of order higher than $ p $.
The remainder term $ \widetilde{\mathcal{R}}_p(s) $ can be expressed using its integral representation
\[
\widetilde{\mathcal{R}}_p(s) = \int_0^s \frac{(s-u)^p}{p!} e^{iHu} \mathrm{ad}_H^{p+1}(O) e^{-iHu} \, du,
\]
with the norm bound 
\[
\|\widetilde{\mathcal{R}}_p(s)\| \leq \int_0^{|s|} \frac{(|s| - u)^p}{p!} \cdot \|e^{iHu}\| \cdot \|\mathrm{ad}_H^{p+1}(O)\| \cdot \|e^{-iHu}\| \, du.
\]

Evaluating the integral gives an upper bound of:
\[
\frac{|s|^{p+1}}{(p+1)!} \cdot \|\mathrm{ad}_H^{p+1}(O)\| = \widetilde{\alpha}_{\text{comm}}\frac{|s|^{p+1}}{(p+1)!}
\]
where
\[\widetilde{\alpha}_{\text{comm}} := \|\mathrm{ad}_H^{p+1}(O)\|.\]
Using this, we have 
\[
\|\widetilde{\mathcal{R}}_p(s)\| = \Or(\widetilde{\alpha}_{\text{comm}}s^{p+1}),
\]
with constants, omitted by the big-O notation, depending only on $p$.

\subsection{Taylor Expansion for the $p$-th Order Trotter Approximation}

To analyze the error associated with the $p$-th order Trotter formula, we now consider the time evolution of the observable $O$ under the Trotter approximation. Recall the approximated time-evolved observable is given by:
\begin{equation}
T_p(\Delta t) := {U_p^{\Delta t}}^\dagger O U_p^{\Delta t}.
\end{equation}
According to the Suzuki formula, $U_p^{\Delta t}$ is expressed as:
\begin{equation}
U_p^{\Delta t} = \REV{e^{-i\Delta t c_1 H_1} e^{-i\Delta t c_2 H_2} \cdots e^{-i\Delta t c_{l-1} H_{l-1}}e^{-i\Delta t c_l H_l}} ,
\end{equation}
where each $ e^{-i\Delta t c_j H_j} $ is an unitary operator, $ H_j \in \{A, B\}$, and the coefficients $ c_j $ are determined by the Suzuki formula for $j \in \{1,\dots,l\}$. The adjoint of this operator is similarly decomposed as:
\begin{equation}
{U_p^{\Delta t}}^\dagger = e^{i\Delta t c_l H_l} e^{i\Delta t c_{l-1} H_{l-1}} \cdots e^{i\Delta t c_2 H_2} e^{i\Delta t c_1 H_1}.
\end{equation}
We now proceed by expanding $ T_p(s) $ in $s$ from the innermost layer of matrix exponential conjugation to the outermost using Taylor's theorem. We will only keep track of those terms of order $O(s^{p+1})$ due to the order condition, and the corresponding commutators. The expansion takes the form:
\begin{equation}
T_p(s) = C_0 + C_1 s + \cdots + C_{p} s^{p} + \mathcal{R}_{p}(s),
\end{equation}
where $ C_0, C_1, \dots, C_p $ are operators independent of $ s $. These terms of order $1,s,\dots, s^p$ will vanish due to the order condition in the final representation of Trotter error. $ \mathcal{R}_{p}(s) $ here denotes the remainder term given by:
\begin{align} \label{eq:remainder_bound}
\mathcal{R}_{p}(s) := & \sum_{k=1}^{l} \sum_{\substack{q_1 + \cdots + q_k = p+1 \\ q_k \neq 0}}
e^{is c_l H_l} \cdots e^{is c_{k+1} H_{k+1}} \nonumber \\
& \quad \cdot \int_0^s \! ds_2 \; e^{is_2 c_k \REV{H}_k} \,
\mathrm{ad}_{{iH_k}}^{\,q_k} \cdots \mathrm{ad}_{{iH_1}}^{\,q_1}(O) \, e^{-is_2 c_k \REV{H}_k} \nonumber \\
& \quad \cdot \frac{(s - s_2)^{q_k-1} \, s^{q_1 + \cdots + q_{k-1}}}{(q_k - 1)! \, q_k! \cdots q_1!} \;
e^{-is c_{k+1} H_{k+1}} \cdots e^{-is c_{l} H_{l}} ,
\end{align}
which can be shown by explicit evaluations of derivatives and iterating the Taylor theorem, as detailed in \cite{ChildsSuTranEtAl2020}. 
Bounding the spectral norm, we obtain
\[
\|\mathcal{R}_{p}(s)\| \leq \sum_{k=1}^{l} \sum_{\substack{q_1 + \dots + q_k = p+1 \\ q_k \neq 0}} 
\left\| \int_0^s ds_2 \, e^{is_2 c_k \REV{H}_k} 
\mathrm{ad}_{H_k}^{q_k} \cdots \mathrm{ad}_{H_1}^{q_1}(O) e^{-is_2 c_k \REV{H}_k} \frac{(s - s_2)^{q_k-1} s^{q_1 + \dots + q_{k-1}}}{(q_k - 1)! q_k! \dots q_1!} \right\|.
\]
After applying the triangle inequality, this gives an upper bound of
\begin{align*}
&\sum_{k=1}^{l} \sum_{\substack{q_1 + \dots + q_k = p+1 \\ q_k \neq 0}} 
\int_0^{|s|} ds_2 \, \frac{(|s| - s_2)^{q_k-1} |s|^{q_1 + \dots + q_{k-1}}}{(q_k - 1)! q_k! \dots q_1!} \cdot 
\|\mathrm{ad}_{H_k}^{q_k} \cdots \mathrm{ad}_{H_1}^{q_1}(O)\| \\
&= \sum_{\substack{q_1 + \dots + q_l = p+1}} 
{ {p+1} \choose {q_1 \cdots q_l} }\frac{|s|^{p+1}}{(p+1)!} \cdot \|\mathrm{ad}_{H_k}^{q_k} \cdots \mathrm{ad}_{H_1}^{q_1}(O)\|. \\
& = \alpha_{\text{comm}}\frac{|s|^{p+1}}{(p+1)!}
\end{align*}
where
\begin{equation}
\alpha_{\text{comm}} := \sum_{\substack{q_1 + \dots + q_k = p+1}} 
{ {p+1} \choose \REV{{q_1 \dots q_k} }}\|\mathrm{ad}_{H_k}^{q_k} \cdots \mathrm{ad}_{H_1}^{q_1}(O)\|.
\end{equation}
Using this, we have
\begin{equation} \label{pth_remainder_bound}
\|\mathcal{R}_{p}(s)\| = \mathcal{O}(\alpha_{\text{comm}}s^{p+1}),
\end{equation}
with the preconstant, omitted by the big-O notation, depending only on $p$. 

\subsection{Commutator Scaling of Observable Errors}\label{subsec:comm_ob}
The local error bound for the $p$-th order Trotter formula is obtained by summing the contributions from the exact remainder term $ \widetilde{\mathcal{R}}_p(\Delta t) $ and the Trotter remainder term $ \mathcal{R}_p(\Delta t) $. From the previous sections, we have:
\begin{equation}
\norm{T_p(\Delta t) - T(\Delta t)} = \norm{\mathcal{R}_p(\Delta t) - \widetilde{\mathcal{R}}_p(\Delta t)}.
\end{equation}
Applying the triangle inequality, we obtain
\begin{equation}
\norm{T_p(\Delta t) - T(\Delta t)} \leq \norm{\mathcal{R}_p(\Delta t)} + \norm{\widetilde{\mathcal{R}}_p(\Delta t)}.
\end{equation}
Using the previously derived bounds in terms of the commutator coefficients, we arrive at the bound
\begin{equation}\label{local_bound_in_commutator}
\norm{T_p(\Delta t) - T(\Delta t)} \leq C_p (\alpha_{\text{comm}} + \widetilde{\alpha}_{\text{comm}}) \Delta t^{p+1},
\end{equation}
where $ C_p $ is a constant depending only on $p$. In particular, terms in $\widetilde{\alpha}_{\text{comm}}$ are upper bounded by terms in $\alpha_{\text{comm}}$ up to a constant by Triangle Inequality:
\[
\|\mathrm{ad}_H^{\,p+1}(O)\|
=\Bigl\|\bigl[A+B,\dots,[A+B,O]\dots\bigr]\Bigr\|
\le\sum_{q=0}^{p+1}\binom{p+1}{q}\|\mathrm{ad}_B^q\mathrm{ad}_A^{\,p+1-q}(O)\|
\le2^{p+1}\,\beta_{\mathrm{comm}},
\]
where $\beta_{\mathrm{comm}}$ is given by \cref{eq:beta}.
The uniformity in $h$ reduces to proving that the commutators are independent of $h$, which is analyzed in detail in the next section and summarized in \cref{thm:nested_poly_obs} for the spectral spatial discretization and \cref{thm:discrete_poly_bound} for the finite difference discretization.

\section{Algebraic Structure and Lemmas}
\label{sec:algebraic_structures}
In this section, we establish the algebraic estimates required for bounding the Trotter observable error for the semiclassical Schr\"odinger operators. For completeness, we begin by introducing the underlying Lie algebra in a spatially continuous setting and reviewing its key properties. We then utilize this algebraic structure to derive the necessary commutator estimates.

While results in the semi-discrete (spatially continuous) setting are generally considered for analyzing spectral discretizations, thanks to their spectral accuracy, this is not the case for finite-difference schemes, which require a fully discrete analysis. Accordingly, we provide detailed estimates in the finite difference setting to validate the approach in a fully discrete context.

\subsection{Height-Width Properties and Lemmas in Semi-discrete Setting}
For our purpose, we consider the set of operators
\[
\mathcal{L}_h = \Bigl\{\sum_{k=0}^n y_k(x)\,h^{m_k}\partial_x^{d_k}\;\Big|\;
n, d_k \in \mathbb{N}, \;  m_k \in \mathbb{Z}, \; k = 0,\dots,n, \; y_k(x) \text{ smooth} \Bigr\}.
\]
We define the height and width for an operator $P \in \mathcal{L}_h$ as follows:
\begin{equation}\label{eq:def_ht_wd}
\mathrm{ht}(P) = \max\{d:h^m \partial_x^d \text{ 
appears in }P\}, \quad \mathrm{wd}(P) = \min\{m: h^m \partial_x^d \text{ appears in }P \}. 
\end{equation}
In particular, we set $\mathrm{ht}(0) = 0$ and $\mathrm{wd}(0) = \infty$.
We have the following structural properties:
\begin{lemma}[Height-Reduction and Width Expansion]\label{lmm:ht_red} 
For any $A, B \in \mathcal{L}_h$ , we have
\begin{equation}
\label{eq:cont_width_addition}
\mathrm{ht}([A,B]) \leq \mathrm{ht}(A) + \mathrm{ht}(B) -1, \quad
\mathrm{wd}([A,B]) \ge \mathrm{wd}(A) + \mathrm{wd}(B).
\end{equation}
\end{lemma}
The proof of \cref{lmm:ht_red} follows from straightforward calculations, which we include in \cref{app:pf_ht_red} for completeness. We note that the height property corresponds to the height reduction in the Lie algebra
\begin{equation}\label{eq:def_Lie_algebrH_l}
\mathcal{L}: = \{\sum_{k=0}^n y_k(x) \partial_x^k, n \in \mathbb{Z}_+, y_0, \cdots, y_n \text{ are smooth}\},
\end{equation}
which contains the Lie algebra generated by $\partial_x^2$ and $V$. This algebra has been carefully studied in the seminal works~\cite{Singh2015,BaderIserlesKropielnickaSingh2014,BaderIserlesKropielnickaSingh2016,IserlesKropielnickaSingh2018,IserlesKropielnickaSingh2019}, leading to a successful new class of Zassenhaus splitting methods and Magnus-based numerical integrators. In particular, the height of an operator can be computed explicitly. For example,
\begin{align*}
[V, \partial_x^2] = -(\partial_x^2 V) - 2(\partial_xV)\partial_x
\end{align*}
is of height 1. The height reduction and width expansion properties play a key role in our proofs.

\begin{lemma}
\label{lmm:single_nested_comm_order} 
Let $C_n=[U_n,[U_{n-1},\dots,[U_2,U_1]\dots]]$
be any grade-$n$ commutator in which each $U_j$ is either $\REV{-\frac{h}{2}} \partial_x^2$ or $h^{-1}V(x)$.  
Let $O_q = y_q(x)h^q\partial_x^q$
be an observable such that $\mathrm{ht}(O) = \mathrm{wd}(O) = q$.
Then
\[
\mathrm{ht}([C_n,O_q])\leq\mathrm{wd}([C_n,O_q]).
\]
\end{lemma}

\begin{proof}
Let $C_n=[U_n,\dots,[U_2,U_1]\dots]$ have $m$ factors of $U_j=h\partial_x^2$ and $n-m$ factors of $U_j=h^{-1}V(x)$.  Then
\[
\sum_{j=1}^n\mathrm{ht}(U_j)
=2m,\qquad
\sum_{j=1}^n\mathrm{wd}(U_j)
= m-(n-m)=2m-n.
\]
By Lemma \ref{lmm:ht_red},
\[
\mathrm{ht}(C_n)\le\sum_{j=1}^n\mathrm{ht}(U_j)-(n-1)
=2m-(n-1),
\qquad
\mathrm{wd}(C_n) \ge \sum_{j=1}^n\mathrm{wd}(U_j)
 = 2m-n.
\]
Since $\mathrm{ht}(O_q)=\mathrm{wd}(O_q)=q$, another application of Lemma \ref{lmm:ht_red} to $[C_n,O_q]$ gives
\[
\mathrm{ht}([C_n,O_q])
\le\mathrm{ht}(C_n)+\mathrm{ht}(O_q)-1
=(2m-n+1)+q-1
=2m-n+q,
\]
while
\[
\mathrm{wd}([C_n,O_q])
\ge \mathrm{wd}(C_n)+\mathrm{wd}(O_q)
=(2m-n)+q.
\]
Hence
\[
\mathrm{ht}([C_n,O_q])\le\mathrm{wd}([C_n,O_q]).
\]
\end{proof}

\begin{theorem}
\label{thm:nested_mono_obs}
Let $C_{n_1}, C_{n_2}, \dots, C_{n_k}$ be any finite sequence of grade-$n_j$ commutators of $-\frac{h}{\REV{2}} \partial_x^2$ and $h^{-1} V$, and let $O_q = y_q(x)h^q\partial_x^q $ be an observable.
Define \[W_k = [C_{n_k},[C_{n_{k-1}}, \dots, [C_{n_1},O_q]\dots]].\]
Then every such nested commutator $W_k$ has that
\[
\mathrm{ht}(W_k)\leq\mathrm{wd}(W_k).
\]
\end{theorem}

\begin{proof}
We will prove by induction on the number $k$ of commutators.

Base case ($k=1$):
\[
W_1 = [C_{n_1},O_q],
\]
and the claim $\mathrm{ht}(W_1)\le\mathrm{wd}(W_1)$ is exactly Lemma \ref{lmm:single_nested_comm_order}.

Inductive step:
Assume the statement holds for $k-1$, so $\mathrm{ht}(W_{k-1})\le\mathrm{wd}(W_{k-1}).$
Set
\[
W_k = [C_{n_k},W_{k-1}].
\]
By Lemma \ref{lmm:ht_red} and the width expansion,
\[
\mathrm{ht}(W_k)\le\mathrm{ht}(C_{n_k})+\mathrm{ht}(W_{k-1})-1,
\quad
\mathrm{wd}(W_k) \ge \mathrm{wd}(C_{n_k})+\mathrm{wd}(W_{k-1}).
\]
Hence
\[
\mathrm{wd}(W_k)-\mathrm{ht}(W_k)
\ge
\bigl[\mathrm{wd}(C_{n_k}) - (\mathrm{ht}(C_{n_k})-1)\bigr]
+ \bigl[\mathrm{wd}(W_{k-1}) - \mathrm{ht}(W_{k-1})\bigr].
\]
By the inductive hypothesis $\mathrm{wd}(W_{k-1}) - \mathrm{ht}(W_{k-1})\ge 0$, so it suffices to show
\[
\mathrm{wd}(C_{n_k}) \ge \mathrm{ht}(C_{n_k})-1.
\]
Write
\[
C_{n_k}
= [U_{n_k},[U_{n_k-1},\dots,[U_2,U_1]\dots]],
\]
with each $U_j$ either $h\partial_x^2$ or $h^{-1}V(x)$.  There are $n_k-1$ nested commutators in this expansion, so applying Lemma \ref{lmm:ht_red} at each of those gives
\[
\mathrm{ht}(C_{n_k})
\le
\sum_{j=1}^{n_k}\mathrm{ht}(U_j)-(n_k-1),
\quad
\mathrm{wd}(C_{n_k})
\ge
\sum_{j=1}^{n_k}\mathrm{wd}(U_j).
\]
Since for each $j$, $\mathrm{ht}(U_j)-\mathrm{wd}(U_j)=1$, we have
\[
\sum_{j=1}^{n_k}\mathrm{ht}(U_j)
=
\sum_{j=1}^{n_k}\mathrm{wd}(U_j)+n_k,
\]
and therefore
\[
\mathrm{ht}(C_{n_k})-1
\le
\Bigl(\sum_{j=1}^{n_k}\mathrm{wd}(U_j)+n_k\Bigr)-n_k
=
\sum_{j=1}^{n_k}\mathrm{wd}(U_j)
\le
\mathrm{wd}(C_{n_k}).
\]

Thus $\mathrm{wd}(C_{n_k})\ge\mathrm{ht}(C_{n_k})-1$, which in turn gives
$\mathrm{wd}(W_k)-\mathrm{ht}(W_k)\ge0$ and equivalently $\mathrm{ht}(W_k)\le\mathrm{wd}(W_k)$, completing the induction.
\end{proof}

The relationship between width and height in the preceding lemmas has a direct implication for the norm of nested commutators after spatial discretization using spectral methods. Specifically in spectral discretization, the spatial derivatives are replaced by the operator $\mathrm{IDFT}_N  \operatorname{diag}((i\xi)^k) \mathrm{DFT}_N$ where $\xi$ denotes the frequency grid, $N$ is the number of spatial grid points, and $\mathrm{IDFT}$ and $\mathrm{DFT}$ are the inverse and forward discrete Fourier transform with $N$ grid points. This representation, along with its exponential, can be efficiently implemented on quantum computers since the DFT and IDFT correspond to the quantum \REV{F}ourier transform (QFT) and inverse QFT, respectively.

As carefully documented in~\cite{BaoJinMarkowich2002,JinMarkowichSparber2011,LasserLubich2020}, to maintain accuracy in the semiclassical regime, the spatial grid size $\Delta x$ for spectral methods must satisfy $\Delta x = \Or(h)$, implying that the number of spatial points scales as $N = \Or(h^{-1})$. Consequently, the spatial discretization of each $k$-th order derivative results in an operator with norm $\Or(h^{-k})$. For instance, the discrete Laplacian $\Delta_x$ has norm $\Or(h^{-2})$, and $\partial_x^k$ has norm $\Or(h^{-k})$.

In other words, a nested commutator term $W$ of height $q$ has norm $\Or(h^{-q})$. If its width is at least $q$, then the structure of the commutators ensures that the norm of $W$ is at most $\Or(1)$.
This leads to the following theorem for spectral spatial discretization:

\begin{theorem}[Uniform-in-$h$ Bound for Nested Commutators of Polynomial Observables]\label{thm:nested_poly_obs}
Let $C_{n_1},\dots,C_{n_k}$ be any sequence of grade-$n_j$ commutators built from $-\frac{h}{\REV{2}}\partial_x^2$ and $h^{-1}V(x)$ satisfying \cref{assump:V}, and let
\[
  O =\sum_{m=0}^q O_m,
  \quad
  O_m = y_m(x)h^m\partial_x^m
\]
be an arbitrary polynomial observable satisfying \cref{assump:observable}.  Define $  W_k^\mathrm{sp} $ as the matrix of
\begin{equation*}
\; [\,C_{n_k},[\,\dots,[\,C_{n_1},O]\,\dots]].
\end{equation*}
after spectral discretizations.
Then $\|W_k^\mathrm{sp}\|_2 = O(h^0)$. In other words,
\[
\|W_k^\mathrm{sp}\|_2 \leq C,
\]
with a constant $C$ depending on $V$ and its derivatives, but independent of $h \in (0,1]$.
\end{theorem}

We finally comment on the case where $O$ is not a finite-degree polynomial in $\partial_x$, but is instead given by an infinite Taylor series $O=\sum_{m\ge 0} y_m(x)h^m\partial_x^m$. Consider $O$ as an $m$-th order mononial $y_m(x)h^m\partial_x^m$, a grade-$n$ nested commutator $\operatorname{ad}_{H_1}\cdots\operatorname{ad}_{H_{n-1}}(O)$, with each $H_j \in \{-\tfrac{h}{2}\partial_x^2,\, h^{-1}V\}$, produces at most $n^m$ growth (including coefficients). For example, the grade-2 commutator $[h^{-1}V(x),\, y(x)h^m\partial_x^m] = -h^{m-1}\sum_{k=1}^m \binom{m}{k}y(x)V^{(k)}(x)\partial_x^{m-k}$ has discretized norm bounded by $\sum_{k=1}^m \binom{m}{k}\le 2^m$. Similarly, the grade-3 commutator $[h^{-1}V,\,[h^{-1}V,\,y(x)h^m\partial_x^m]] = y(x)h^{m-2}\sum_{j=1}^m\sum_{\ell=1}^{m-j}\binom{m}{j}\binom{m-j}{\ell}V^{(j)}V^{(\ell)}\partial_x^{m-j-\ell}$ has discretized norm bounded by $3^m$; other grade-3 patterns grow more slowly. In general, the largest norm of a grade-$n$ commutator is bounded by $n^m$, and there are at most $2^{n-1}$ possible choices of grade-$n$ nested commutators. Hence the total grade-$n$ contribution is bounded by $2^{n-1}n^m$. 
For a non-polynomial observable $O$ defined by an infinite Taylor series, each individual term after nested commutators still remains uniformly bounded in $h$ by the theorem above. To keep the entire infinite series uniformly bounded, it suffices that the Taylor coefficients, after multiplication by $2^{n-1}n^m$, remain summable. This is satisfied, for example, when the Taylor coefficients exhibit factorial decay (as in the Taylor coefficients of an exponential function). Hence non-polynomial observables with sufficiently fast decaying Taylor coefficients may still satisfy uniform-in-$h$ bounds.

\subsection{Spatial Discretization and Lemmas in Fully Discrete Setting}
In this section, we analyze the algebraic structure in the fully discrete setting using finite difference spatial discretization. Compared to the spectral methods discussed in \cref{thm:nested_poly_obs}, the finite difference approach requires additional care, as standard properties such as the product rule do not hold in the same form as they do at the level of continuous operators.

For example, consider the forward difference operator $\Delta_+$ defined by
\begin{equation}
(\Delta_+ u)_\ell = \frac{u_{\ell+1}-u_\ell}{\Delta x}, 
\end{equation}
with spatial grid size $\Delta x$.
The corresponding discrete product rule takes the form
\begin{equation}
    \Delta_+(uv)_\ell =  \left(\Delta_+ u\right)_\ell \, v_{\ell+1} +  u_\ell \, \left(\Delta_+ v\right)_\ell .
\end{equation}
where the index shift in the first term differs from the product rule in the continuous setting.

We work on a uniform grid of $N$ points
\[
  x_j = a + (b-a)\,\frac{j}{N},\quad j=0,1,\dots,N-1,
\]
with periodic boundary conditions.
Define the scaled finite-difference matrices
\[
D_F = N
\begin{bmatrix}
-1 & 1 & 0 & \cdots & 0\\
0 & -1 & 1 & \cdots & 0\\
\vdots & & \ddots & \ddots & 1\\
1 & 0 & \cdots & 0 & -1
\end{bmatrix},
\quad
D_B = -D_F^\dagger,
\quad
D_2 = N^2
\begin{bmatrix}
-2 & 1 & 0 & \cdots & 1\\
1 & -2 & 1 & \cdots & 0\\
\vdots & & \ddots & \ddots & 1\\
1 & 0 & \cdots & 1 & -2
\end{bmatrix}.
\]
One checks easily that
\[
D_2 = D_B\,D_F = D_F\,D_B,
\]
and that $D_F$ and $D_B$ are normal and commute.  We now define higher‐order differences by
\begin{equation}
\label{eq:Dk_definition}
D_k \;=\;
\begin{cases}
D_2^{\,k/2},&k\text{ even},\\
D_F\,D_2^{\,(k-1)/2},&k\text{ odd}.
\end{cases}
\end{equation}
Since $D_F,D_B,D_2$ all commute, one sees immediately that
\[
[D_k,D_j]=0 \quad\forall j,k\ge1.
\]
We only show result with this definition in the following proofs, but the same follows if we define odd powers $D_k$ with $D_BD^{\frac{k-1}{2}}$ without the loss of generality.

Let $y(x)$ be any smooth periodic function and set
\[
Y = \mathrm{diag}\bigl(y(x_0),y(x_1),\dots,y(x_{N-1})\bigr).
\]
Similarly, $Y_k$ can be defined according to the function $y_k(x)$.
Clearly $Y_k$ commutes with $Y_j$ for any such matrices. 

We say an $N\times N$ matrix $P$ has \emph{height} $m$ if $m$ is the highest order such that $||P||_2$ grows as $N^m$ as $N \to \infty$. By the definition of $D_F$, $D_B$, and $D_k$, we have
\[
\mathrm{ht}(D_F) = 1, \quad \mathrm{ht}(D_B) = 1, \quad \mathrm{ht}(D_k) = k.
\]
We will now show commuting a single difference matrix with a diagonal multiplier reduces height by at least 1. 

\begin{lemma}[Commutator of $D_k$ with $Y$]
\label{lmm:discrete_Dk_with_y}
For $k\ge1$,
\[
[D_k,Y]
=D_kY - YD_k
\]
has height at most $k-1$, i.e.
\[
\mathrm{ht}([D_k,Y])\leq k-1.
\]
\end{lemma}

\begin{proof}
For base case $k=1$,  by Mean Value Theorem there is some $\xi_j$ between $x_j$ and $x_{j+1}$ so that
\[
  y(x_{j+1}) - y(x_j)
  = y'(\xi_j)\,(x_{j+1}-x_j)
  = \frac{y'(\xi_j)}{N}.
\]
Thus for any vector $u$,
\[
  ([D_F,Y]\,u)_j
  = N\bigl(y(x_{j+1})-y(x_j)\bigr)\,u_{j+1}
  = y'(\xi_j)\,u_{j+1},
\]
and so
\[
  \|[D_F,Y]\|_2 \;\le\;\max_x|y'(x)|,
\] which means $\mathrm{ht}([D_F,Y])\le0$.

For the inductive step, suppose $\mathrm{ht}([D_{k-1},Y])\le k-2$.  Then we have the following cases:

\emph{Case 1: $k$ is even.}  Write $k=2m$.  Then by definition
\[
D_k \;=\; D_2^m \;=\; D_B\,D_{k-1}.
\]
Hence
\[
[D_k,Y]
= [\,D_B\,D_{k-1},\,Y]
= D_B\,[D_{k-1},Y] \;+\; [D_B,Y]\,D_{k-1}.
\]
Since
\[
\mathrm{ht}(D_B)=1,\quad
\mathrm{ht}\bigl([D_{k-1},Y]\bigr)\le k-2,\quad
\mathrm{ht}([D_B,Y])\le0,\quad
\mathrm{ht}(D_{k-1})= k-1,
\]
we obtain
\[
\mathrm{ht}\bigl(D_B\,[D_{k-1},Y]\bigr)\le 1+(k-2)=k-1,
\quad
\mathrm{ht}\bigl([D_B,Y]\,D_{k-1}\bigr)\le 0+(k-1)=k-1.
\]
Therefore $\mathrm{ht}([D_k,Y])\le k-1$ in the even case.

\medskip
\emph{Case 2: $k$ is odd.}  Write $k=2m+1$.  Then
\[
D_k \;=\; D_F\,D_2^m \;=\; D_F\,D_{k-1},
\]
and
\[
[D_k,Y]
= D_F\,[D_{k-1},Y] \;+\; [D_F,Y]\,D_{k-1}.
\]
Using
\[
\mathrm{ht}(D_F)=1,\quad
\mathrm{ht}\bigl([D_{k-1},Y]\bigr)\le k-2,\quad
\mathrm{ht}([D_F,Y])\le0,\quad
\mathrm{ht}(D_{k-1})= k-1,
\]
we similarly get
\[
\mathrm{ht}\bigl(D_F\,[D_{k-1},Y]\bigr)\le 1+(k-2)=k-1,
\quad
\mathrm{ht}\bigl([D_F,Y]\,D_{k-1}\bigr)\le 0+(k-1)=k-1,
\]
and hence $\mathrm{ht}([D_k,Y])\le k-1$. This finishes the induction.
\end{proof}

\begin{lemma}
\label{lmm:discrete_ht_reduction}
Let $k,j\ge1$ and let $Y_k,Y_j$ be any two smooth diagonal matrices.  Then
\[
\mathrm{ht}\bigl([Y_kD_k,Y_jD_j]\bigr)\le\;k+j-1.
\]
\end{lemma}

\begin{proof}
Using commutator identities, we simplify
\begin{align*}
[Y_k D_k,Y_j D_j] = & [Y_k D_k,Y_j]D_j + Y_j[Y_k D_k, D_j] 
\\
= & Y_k[ D_k,Y_j]D_j + [Y_k ,Y_j]D_k D_j + Y_j Y_k [ D_k, D_j] + Y_j[Y_k, D_j] D_k. 
\end{align*}
Since $[ D_k, D_j] =  [Y_k ,Y_j] = 0 ,$ this simplifies to 
\[
[Y_k D_k,Y_j D_j]
= Y_k[ D_k,Y_j]D_j + Y_j[Y_k, D_j] D_k. 
\]
So it remains to show that $Y_k[ D_k,Y_j]D_j + Y_j[Y_k, D_j] D_k$ doesn't contain a height $k+j$ term. 
 
By Lemma \ref{lmm:discrete_Dk_with_y}, $\mathrm{ht}([D_k,Y_j])\le k-1$ and $\mathrm{ht}([Y_k,D_j])\le j-1$, while $\mathrm{ht}(D_j) = j$, $\mathrm{ht}(D_k) = k$, and $\mathrm{ht}(Y_\ell)=0$. Therefore each of the two remaining terms has 
\[
\mathrm{ht}\bigl(Y_k[ D_k,Y_j]D_j\bigr) \le k+j-1, \quad \mathrm{ht}\bigl(Y_j[Y_k, D_j] D_k\bigr) \le k+j-1.
\]
and the lemma follows.
\end{proof}

Consider the set of operators
\[
\widetilde{\mathcal{L}}_h = \Bigl\{\sum_{k=0}^n Y_k\,h^{m_k}D_{d_k}\;\Big|\;
n, d_k \in \mathbb{N}, \;  m_k \in \mathbb{Z}, \; k = 0,\dots,n \Bigr\}.
\]
We define, for $P \in \widetilde{\mathcal{L}}_h$ that
\[
\mathrm{ht}(P)=\max_k\{d_k: h^{m_k}D_{d_k}\text{ appears in }P\},
\quad
\mathrm{wd}(P)=\min_k\{m_k: h^{m_k}D_{d_k}\text{ appears in }P\}.
\]
We set $\mathrm{ht}(0) = 0$ and $\mathrm{wd}(0) = \infty$. This height definition is consistent with the previous notion based on the growth of $N$, as the norm of $D_{d_k}$ grows as $N^{d_k}$ as $N \to \infty$. 
\begin{lemma}[Discrete Height–Reduction and Width–Expansion]
\label{lmm:discrete_ht_and_width}
For any $P,Q\in\widetilde{\mathcal{L}}_h$,
\[
\mathrm{ht}([P,Q])\le \mathrm{ht}(P)+\mathrm{ht}(Q)-1,
\qquad
\mathrm{wd}([P,Q]) \geq \mathrm{wd}(P)+\mathrm{wd}(Q).
\]
\end{lemma}
\begin{proof}
This is immediate by Lemma \ref{lmm:discrete_ht_reduction} and multiplication of lowest power of $h$.
\end{proof}

\begin{lemma}
\label{lmm:single_comm_discrete}
Let
$C_n = [\,U_n,[\,U_{n-1},\dots,[\,U_2,U_1]\dots]\,]$
be any nested commutator of grade $n$, where each $U_j$ is either
\begin{equation}\label{eq:fd_A_B}
A = -\frac{h}{\REV{2}}D_2,
\quad
B = \frac{1}{h}\mathrm{diag}\bigl(V(x_j)\bigr).
\end{equation}
Let $O_q= Y_q h^qD_q$ be an observable with $\mathrm{ht}(O_q)=\mathrm{wd}(O_q)=q$. Then
\[
\mathrm{ht}([C_n,O_q]) \le \mathrm{wd}([C_n,O_q]).
\]
\end{lemma}

\begin{proof}
Write $C_n$ with $m$ copies of $A$ and $n-m$ of $B$.  Since
$\mathrm{ht}(A) = 2$, $\mathrm{wd}(A)=1$ and
$\mathrm{ht}(B) = 0$, $\mathrm{wd}(B)=-1$,  \cref{lmm:discrete_ht_and_width}  gives
\begin{equation}
\label{eq:discrete_ht_wd_Cn}
\mathrm{ht}(C_n) \le 2m-(n-1),
\quad
\mathrm{wd}(C_n)\ge2m-n
\end{equation}
for any such $C_n.$
Applying \cref{lmm:discrete_ht_and_width} again to $[C_n,O_q]$ (with $\mathrm{ht}(O_q) = q, \, \mathrm{wd}(O_q)=q$) yields
\[
\mathrm{ht}([C_n,O_q])\le(2m-(n-1))+q-1=2m-n+q,
\quad
\mathrm{wd}([C_n,O_q])\ge 2m-n+q,
\]
and hence $\mathrm{ht}([C_n,O_q])\le\mathrm{wd}([C_n,O_q])$.

\end{proof}

\begin{theorem}
\label{thm:discrete_nested_mono_obs}
Let $C_{n_1},\dots,C_{n_k}$ be any finite sequence of grade-$n_j$ commutators of $A$ and $B$, and let
$O_q= Y_q h^qD_q$
be an observable.  Define
$W_k = [\,C_{n_k},[\,C_{n_{k-1}},\dots,[\,C_{n_1},O_q]\dots]\,].$
Then every such nested commutator satisfies
\[
\mathrm{ht}(W_k)\le \mathrm{wd}(W_k).
\]
\end{theorem}

\begin{proof}
The base case $k=1$ is  \cref{lmm:single_comm_discrete}.  Assume the result holds for $k-1$, so
$\mathrm{ht}(W_{k-1})\le\mathrm{wd}(W_{k-1})$.  Set
$W_k=[C_{n_k},W_{k-1}]$.  By \cref{lmm:discrete_ht_and_width},
\[
\mathrm{ht}(W_k)\le\mathrm{ht}(C_{n_k})+\mathrm{ht}(W_{k-1})-1,
\quad
\mathrm{wd}(W_k)\ge \mathrm{wd}(C_{n_k})+\mathrm{wd}(W_{k-1}).
\]
Since $\mathrm{ht}(C_{n_k})-1 \le   \mathrm{wd}(C_{n_k})$ by \cref{eq:discrete_ht_wd_Cn} , it follows that
$\mathrm{ht}(W_k)\le \mathrm{wd}(W_k)$, completing the induction.
\end{proof}

\begin{theorem}[Uniform-in-$h$ Bound for Polynomial Observables in the Fully Discrete Setting]\label{thm:discrete_poly_bound}
Let $C_{n_1},\dots,C_{n_k}$ be any finite sequence of grade-$n_j$ commutators of $A$ and $B$ \REV{given by \cref{eq:fd_A_B} with the potential $V$ satisfying \cref{assump:V}}, and let
\[
  O =\sum_{m=0}^q O_m,
  \quad
  O_m= Y_m\,h^m\,D_m 
\]
be a polynomial observable \REV{whose coefficients satisfies \cref{assump:observable}}.  Define the nested commutator
\[
  W_k \;=\; [\,C_{n_k},[\,C_{n_{k-1}},\dots,[\,C_{n_1},O]\,\dots]\,].
\]
Then
$  \|W_k\|_2 = O(1)$
uniformly in $h \in (0,1]$. In other words,
\[
 \|W_k\|_2 \leq C, 
\]
with a constant $C$ depending on $V$ and its derivatives, but independent of $h \in (0,1]$.
\end{theorem}
\begin{proof}
By linearity of the commutator,
\[
  W_k
  = \sum_{m=0}^q [\,C_{n_k},[\,\dots,[C_{n_1},O_m]\dots]\,].
\]
But each summand is exactly the nested commutator of the monomial $O_m$, so by Theorem \ref{thm:discrete_nested_mono_obs}
$\mathrm{ht}\le\mathrm{wd}$ for that term and hence its operator norm is $\REV{\Or}(1)$.  A single application of the triangle inequality then gives
\[
  \|W_k\|_2
  \;\le\;
  \sum_{m=0}^q
    \bigl\|\,[C_{n_k},[\,\dots,[C_{n_1},O_m]\dots]]\bigr\|_2
  = \REV{\Or}(1),
\]
as claimed.
\end{proof}

\begin{remark}
In the fully discrete setting this shows that any finite polynomial observable
 $\;O=\sum Y_m h^m D_m$\,
yields nested commutators whose norms remain uniformly bounded in the small-$h$ limit, exactly as in the continuous case.
\end{remark}

\section{Numerical Comparison}
\label{sec:numerics}
While our theoretical results provide rigorous error bounds for higher-order Trotter formulas, we present numerical simulations to further validate these findings. The primary purpose of this section is to illustrate the behavior of higher-order Trotter-Suzuki formulas together with spatial discretization. We consider the semiclassical Schrödinger equation for the Hamiltonian:
\begin{equation}
    -\frac{h}{2} \Delta + \frac{1}{h}V(x), \quad V(x) = \cos(x), \quad x \in [-\pi, \pi],
\end{equation}
with periodic boundary conditions. The operators $A, B, H$ and $O$ in this session follow the finite difference discretization. For example, the Laplacian $ -\Delta $ is discretized using a second-order finite difference scheme, and the potential $ V(x) $ is represented as a diagonal matrix.

\begin{figure}[htbp]
    \centering
    \subfloat[Observable Errors vs. Step Size $\Delta t$]{
        \includegraphics[width=.45\textwidth]{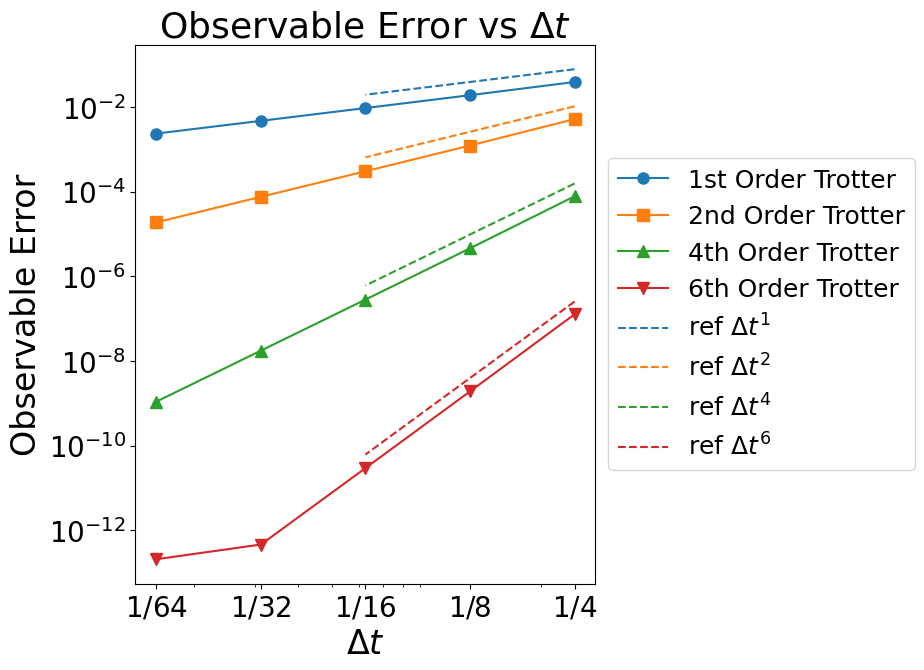}
        \label{fig:observable}
    }
    \subfloat[Unitary Errors vs. Step Size $\Delta t$]{
        \includegraphics[width=.45\textwidth]{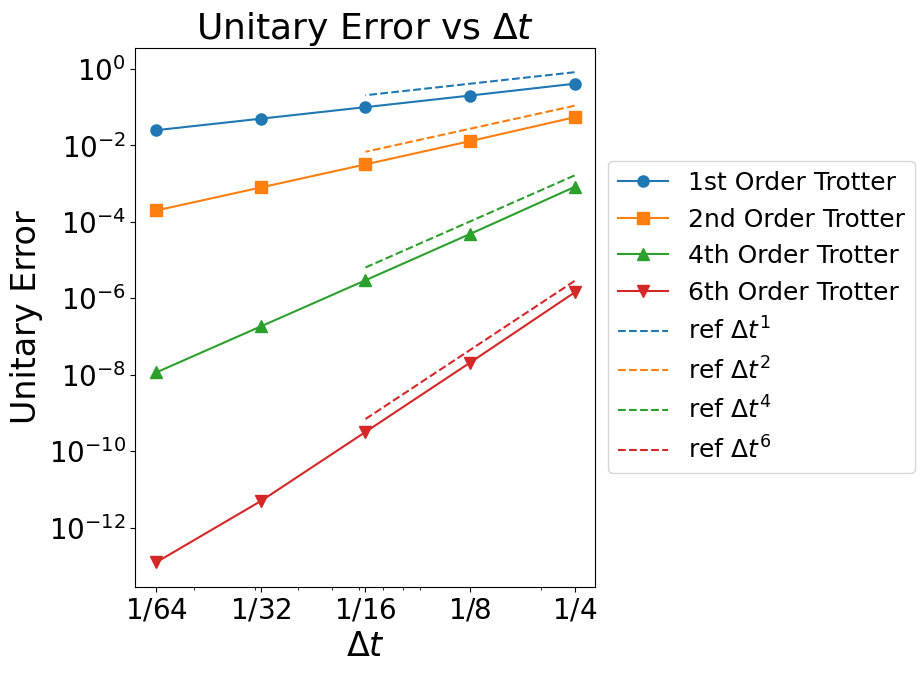}
        \label{fig:unitary}
    }
    \caption{Log-log plots showing the convergence of errors with respect to the step size $\Delta t$ for 1st, 2nd, 4th, and 6th-order Trotter methods. (a) Observable errors are plotted against $\Delta t$, demonstrating convergence consistent with the expected asymptotic scaling. (b) Unitary errors are plotted against $\Delta t$, also showing the convergence rate aligning with the order of the Trotter method. The $\Delta t^p$ reference lines provide a visual confirmation of the expected scaling behavior for each method.}
    \label{fig:errors_vs_dt}
\end{figure}

We first examine the convergence rate with respect to the time step $ \Delta t $ for 1st, 2nd, 4th, and 6th-order Trotter formulas. The time steps are chosen as $ \Delta t = 1/4, 1/8, 1/16, 1/32, 1/64 $, and the simulations run until the final time $ t = 0.5 $. The semiclassical parameter is set to $ h = 1/64 $, and the spatial grid uses $ N = 64 $ points. Figure~\ref{fig:errors_vs_dt} presents the results for both observable and unitary errors as a function of $ \Delta t $. We observe that the observable and unitary errors converge at rates consistent with the order of the Trotter method, with reference lines $ \Delta t^p $ providing visual confirmation of the expected scaling behavior for each method as shown in \cref{thm:global_error}.

\begin{figure}[htbp]
    \centering
    {
        \includegraphics[width=.3\textwidth]{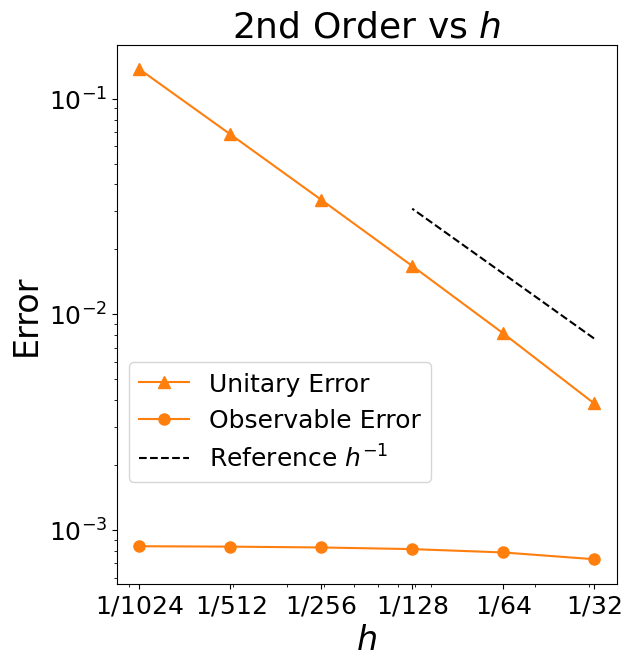}
        \label{fig:2nd_order_h}
    }
    {
        \includegraphics[width=.3\textwidth]{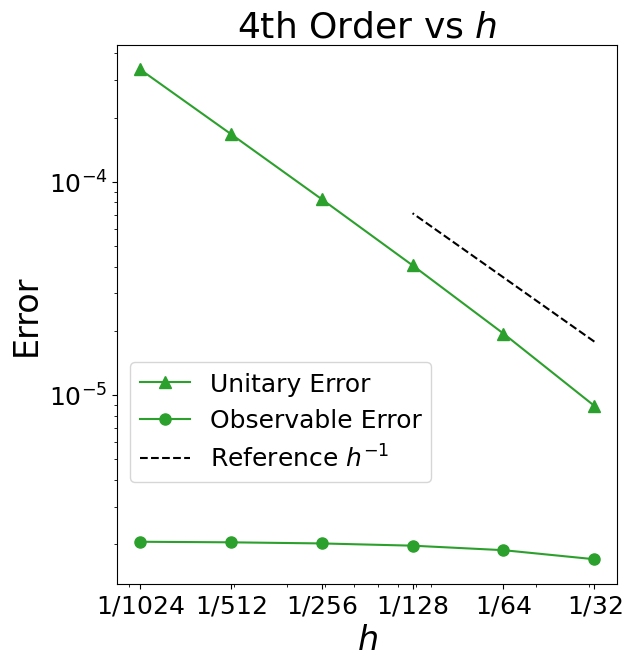}
        \label{fig:4th_order_h}
    }
    {
        \includegraphics[width=.3\textwidth]{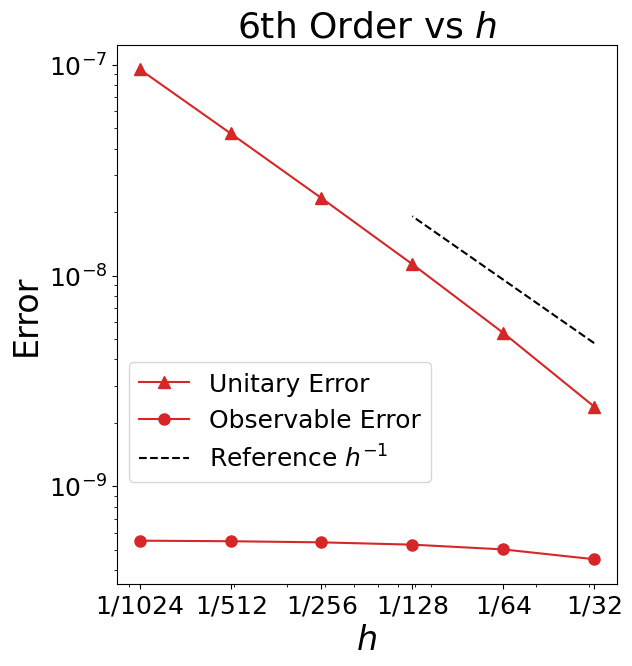}
        \label{fig:6th_order_h}
    }

    \caption{Log-log plots showing unitary and observable errors as a function of the semiclassical parameter $ h $ for 2nd, 4th, and 6th-order Trotter methods. Each plot demonstrates that while the unitary error scales as $ h^{-1} $ (indicated by the reference $ h^{-1} $ line), the observable error remains independent of $ h $, with no apparent slope. This highlights the observable error’s insensitivity to the semiclassical parameter $ h $, particularly scaling independently with $h^{-1}$ aligning with our theoretical results.}
    \label{fig:errors_vs_h}
\end{figure}

Next, we investigate how the errors depend on the semiclassical parameter $ h $. Keeping the time step fixed at $ \Delta t = 0.1 $, we vary $ h $ over $ h = 1/32, 1/64, 1/128, 1/256, 1/512, 1/1024 $. Figure~\ref{fig:errors_vs_h} shows the unitary and observable errors versus $ h $ for the Trotter method. The unitary error increases as $ h $ decreases, approximately scaling with $ h^{-1} $, while the observable error remains relatively constant. This indicates that the observable error is independent of $ h^{-1} $, aligning with \cref{thm:global_error}.

\begin{figure}[H]
    \centering
    {\includegraphics[width=.55\textwidth]{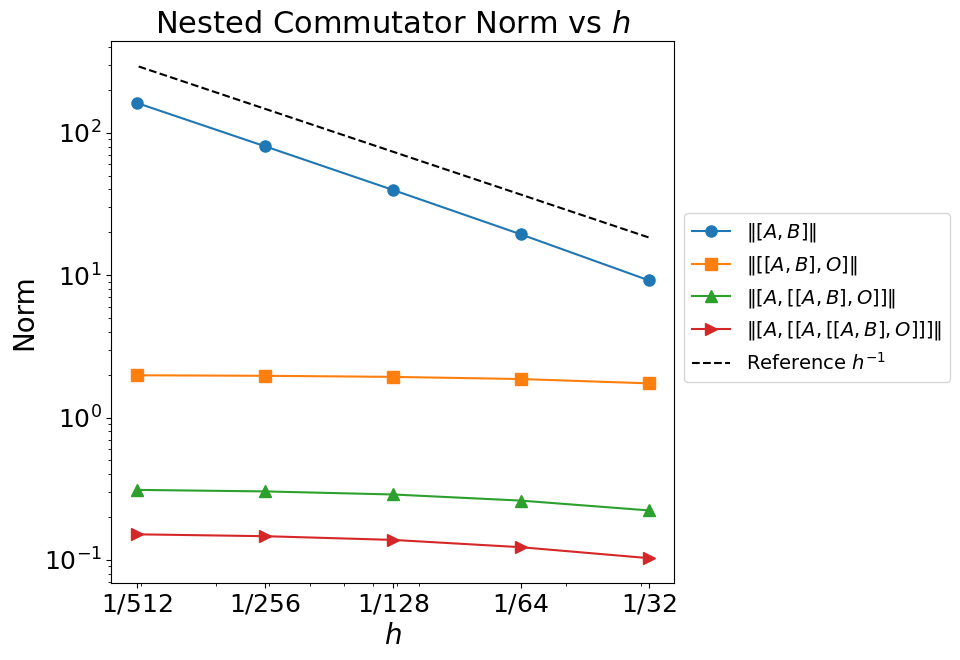}}
    \caption{Log-log plot of commutators illustrating the scaling behavior of nested commutators with observables. The norm of the first commutator $[A, B]$ scales as 
    $\REV{\Or}(1/h)$, confirming the dependence on the semiclassical parameter $h^{-1}$. In contrast, higher-order nested commutators with the observable $O$, such as $[[A, B], O]$, $[A, [[A, B], O]]$, and $[A, [A, [[A, B], O]]]$, demonstrate independence from $h^{-1}$ and remain consistent with $\REV{\Or}(1)$ scaling. The observed behavior aligns with theoretical predictions of observable error insensitivity to $h^{-1}$ in higher-order Trotter approximations.}
    \label{fig:commutators}
\end{figure}

Finally, we observe the behavior of commutators and nested commutators involving the observable $ O $. The commutator norm $\|[A, B]\|$,
scales as $\REV{\Or}(1/h) $, while the nested commutators $[[A, B], O]$, $[A, [[A, B], O]]$, and higher orders are expected to scale as $\REV{\Or}(1) $, independent of $ h^{-1} $.
Figure~\ref{fig:commutators} presents the log-log plots of the norms of $[A, B]$, $[[A, B], O]$, $[A, [[A, B], O]]$, and $[A, [A, [[A, B], O]]]$ as functions of $ h $. The numerical results confirm two key findings consistent with \cref{thm:discrete_poly_bound}. First, the norm of the commutator $[A, B]$ scales as $\Or(1/h)$, in agreement with theoretical predictions. Second, nested commutators involving the observable $O$ remain bounded, exhibiting $\Or(1)$ scaling and thus demonstrating independence from the $h^{-1}$ scaling.

\section{Conclusion and Discussion}
\label{sec:conclusion}

In this work, we present the first purely algebraic proof of uniform-in-\(h\) observable error bounds for arbitrarily high-order Trotter–Suzuki formulas applied to the semiclassical Schrödinger equation. By directly comparing the exact and Trotterized evolution of a polynomial observable \(O\), we avoid any reliance on Egorov-type theorems or semiclassical limits, and thereby retain the full $p$th-order convergence in $\Delta t$ without introducing additive $\REV{\Or}(h)$ terms.  Our main result (\cref{thm:global_error}) shows that 
\[
\bigl\|\, (U_p^{\Delta t})^{\dagger n}O\,(U_p^{\Delta t})^n \,-\, e^{iHt}O\,e^{-iHt}\bigr\|\;\le\;C\,\Delta t^p,
\]
with $C$ independent of $h^{-1}$.  As a byproduct, \cref{prop:general_observable_error} yields a general nested-commutator bound valid for any decomposition $H=A+B$ and bounded $O$ beyond the semiclassical regime.

Key to our analysis is a height-width property in the Lie algebra generated by the kinetic and potential parts, which we first establish in the spectral (semi-discrete) setting and then extend to the finite-difference discretization.  These algebraic estimates guarantee that all $(p+1)$-fold commutators appearing in the local Trotter error remain uniformly bounded as $h\to0$.  Numerical experiments in \cref{sec:numerics} confirm the predicted $\Delta t^p$ rates, the $h$-independence of the observable error, and the boundedness of the nested commutators.

An immediate future direction is to extend our algebraic framework to time‐dependent Hamiltonian simulation, analyzing observable‐error scaling under generalized Trotter formulas~\cite{HuyghebaertDeRaedt1990,AnFangLin2021,MizutaTatsuhikoKeisuke2024} and exploring their applications in driven or adiabatic quantum algorithms. The algebraic structure investigation could also potentially be helpful for improving the complexity estimate of quantum Magnus algorithms, due to their cost dependence on time-dependent nested commutators~\cite{AnFangLin2022,FangLiuSarkar2025,BornsweilFangZhang2025}, which will be left as a future direction.

\REV{In this work, we consider potentials in $S(1)$ (see \cref{assump:V}). The (one-body or many-body) Coulomb potentials, for instance, fall outside this setting. Interestingly, the Trotter convergence for Coulomb interactions is substantially more delicate: one can obtain at most a one-quarter convergence rate for the first-order Trotter formula, and this rate is sharp, being saturated by the ground state of the hydrogen atom as the initial condition~\cite{FangWuSoffer2025,BurgarthFacchiHahnJohnssonYuasa2024,BeckerGalkeSalzmannLuijk2024}. For such singular potentials, higher-order Trotter–Suzuki splittings may behave differently depending on the order in which $A$ and $B$ are applied. This contrasts with our setting, where the regularity of $V$ ensures that the ordering of $A$ and $B$ does not matter (in terms of the uniformity of $h$). Extending the uniform-in-$h$ observable result to Coulomb potentials is an interesting direction for future work.}

\section*{Acknowledgements}
The authors thank Francois Golse for valuable comments. This work is supported by National Science Foundation via the grant DMS-2347791 and DMS-2438074 (D.F. and C.Q.), and the U.S. Department of Energy, Office of Science, Accelerated Research in Quantum Computing Centers, Quantum Utility through Advanced Computational Quantum Algorithms, grant no. DE-SC0025572 (D.F.). 

\appendix

\section{Proof of \cref{lmm:ht_red}}\label{app:pf_ht_red}

\begin{proof}
    Let $A, B \in \mathcal{L}_h$. For the height reduction, it suffices to show the highest order $\partial_x$ term cancels out in the commutator $[A,B]$.
Let any term in $A$ with the derivative order equal to its height be
\begin{equation}
    y(x) h^m \partial_x^d, 
\end{equation}
and similarly let that in $B$ be
\begin{equation}
    z(x) h^k \partial_x^j.
\end{equation}
Then their commutator is given by
\begin{equation}
    [y(x) h^m \partial_x^d, z(x)  h^k \partial_x^j]
    = h^{m+k} \left( y(x)  [ \partial_x^d, z(x)]\partial_x^j
    +  z(x)  [y(x) , \partial_x^j] \partial_x^d \right).
\end{equation}
By the product rule, we have
\begin{equation}
    [ \partial_x^d, z(x)] = z(x)\partial_x^d  + \left(\text{terms of height } \leq d-1\right) - z(x)\partial_x^d = \left(\text{terms of height } \leq d-1 \right),  
\end{equation}
and a similar reduction holds for $[y(x) , \partial_x^j]$.
This shows that the highest-order term in the derivative (involving $\partial_x^{d+j}$) cancels out in the commutator.

The width expansion follows from the multiplication of the $h$-powers in $A$ and $B$. When $[A, B] = 0$, the property holds trivially since $\mathrm{wd}(0) = \infty$. Otherwise, the leading terms (i.e., those with the smallest width) in both $A$ and $B$ either do not commute -- in which case the width of the resulting commutator equals the sum of the widths of $A$ and $B$ -- or commute, in which case the width increases beyond the sum.
\end{proof}

\bibliographystyle{unsrtnat}
\bibliography{ref}

\end{document}